\documentclass{amsart}
\usepackage{amssymb}

\begin{document}

\bibliographystyle{plain}

\title[ R.F. Shamoyan, N.M. Makhina]{Sharp extensions of analytic function spaces in the unit disk and related domains}

\author[]{R.F. Shamoyan}
\address{Bryansk State University, Bryansk, Russia}
\email{\rm rshamoyan@gmail.com}

\author[]{N.M. Makhina}

\address{Bryansk State University, Bryansk, Russia}
\email{\rm makhina32@gmail.com}

\date{}

\maketitle

\begin{abstract}
In this expository paper we collect many recent advances in analytic function spaces of several complex variables related with trace problem. We consider various function space of analytic functions of several variables in various domains in ${\rm {\mathbb C}}^{n} $ and provide or complete descriptions of traces or estimates of traces of various analytic function spaces in various domains obtained in recent years, by various authors. The problem to find sharp estimates of Hardy analytic function spaces in the unit polydisk first was posed by W. Rudin in 1969. Since then many papers appeared in literature. We collect in this expository paper not only already many known results on traces of various analytic function spaces in product domains but also discuss various new interesting results related with this problem. Related to trace problem various results were provided previously by G. Henkin, E. Amar, H. Alexander and various other authors. Finally note that our trace theorems are closely related with the Bergman type projections acting between function spaces with different dimensions.
\end{abstract}

\footnotetext[1]{Mathematics Subject Classification 2010 Primary 32A10, Secondary 32T15.  Key words
and Phrases: Trace operator, Bergman and Hardy mixed norm spaces, Herz spaces, Bergman type projection, polydisk, polyball, tube domains, pseudoconvex domains, analytic functions.}

\section{INTRODUCTION}

Let ${\rm {\mathbb C}}$ denote the set of complex numbers and let ${\rm {\mathbb C}}^{n} ={\rm {\mathbb C}}\times ...\times {\rm {\mathbb C}}$ denote the Euclidean space of complex dimension n. Let $U$ be the unit disk on a complex plane ${\rm {\mathbb C}}$ and $U^{n} =U\times ...\times U$ be the unit polydisk a bounded domain in ${\rm {\mathbb C}}^{n} $.

Let $f$ be an analytic function in the unit polydisk from a certain fixed quazinormed analytic function space $X(U^{n} )$ in the unit polydisk. This expository paper concerns the following rather simple problem to find precise estimates of $|f(z,...,z)|$ set where $z$ is an arbitrary point from a unit disk in ${\rm {\mathbb C}}^{n} $.

Obviously, very similarly such a problem may be posed in tubular domains over symmetric cones and bounded strongly pseudoconvex domains with smooth boundary and in products of such type domains in ${\rm {\mathbb C}}^{n} $.

This problem was posed for the first time in W. Rudins book [1] in 1969 for Hardy spaces in bidisk and since then many papers appeared in this direction. We mention for example old papers [2-9] and various references there.

First the problem was posed for classical Hardy spaces in the unit polydisk then various authors studied this type problems in various other analytic function spaces in products of unit ball, in polydisk, tubular domains over symmetric cones, and bounded strongly pseudoconvex domains with smooth boundary.

Various new interesting estimates, and sharp results were obtained in this direction in recent years. See for example recent papers [10-40] and also many references there also.

The goal of this paper to collect all these new interesting results in this new expository research paper. Note that various new open interesting problems related with this problem will also be provided in this paper.

We start with Trace results first in the simplest domain namely the unit polydisk then provide new results in more complicated domains in ${\rm {\mathbb C}}^{n} $ such as tube domains over symmetric cones and bounded strongly pseudoconvex domains with smooth boundary in ${\rm {\mathbb C}}^{n} $.

We need first the following precise definition of Traces in analytic function spaces in the unit polydisk in ${\rm {\mathbb C}}^{n} $. Note we provide a definition in the unit polydisk but very similarly we can pose such a problem in other more complicated domains in ${\rm {\mathbb C}}^{n} $. We leave this easy task to readers.

\textbf{Definition 1.} \textit{Let $X(U^{n} )$ be a certain fixed quazinormed function space in the unit polydisk in ${\rm {\mathbb C}}^{n} $. We say that $Trace(X)=Y$, where $Y$ is an analytic function space in the unit disk in ${\mathbb C},$ if for every $f$ function from $X$ space $f(z,...,z)$ function belongs to $Y$ space of analytic functions in the unit disk and the reverse is also true, for every $f$ from $Y$ space we can find an $F$ function from $X$ analytic function space so that we have $F(z,...,z)=f(z)$ for all $z$ points from a unit disk $U$on a complex plane ${\rm {\mathbb C}}$. }

Note, in very general tubular domains and bounded strongly pseudoconvex domains we can define traces of analytic functions very similarly as in the polydisk which we provided in our definition below. We refer for basic facts on these domains and analytic spaces on them to [13, 20, 30-34, 41-44]. We shortly pose in this expository paper also some open problems related with traces of analytic function spaces in these domains.

We note that we replace $TraceY$ below by simpler $DY$ sometimes.

In this expository paper we will collect many old and new results concerning this problem for many \textit{X }analytic function spaces in various product domains in ${\rm {\mathbb C}}^{n} $ (not only the unit polydisk). First below we provide definitions needed for formulations of such type results then we formulate our results (old and new) on Traces.

\textbf{Remark 1.} We also note that very similar problems Trace theorems or various estimates of traces of functions were given, proved and considered in many complicated functional spaces in ${\rm {\mathbb R}}^{n} $ by many authors. We refer the reader for example to [45-47] and many references there. 

Trace theorems have also interesting extensions and also nice applications, we for example refer the reader to papers of E. Amar, C. Menini [48] and D. Clark [49] where such type applications and such type extensions were provided.

Finally note that our trace theorems are closely related with the Bergman type projections acting between function spaces with different dimensions.

For the convenience of the reader we organize this expository paper in the following simple manner. We divide all results concerning traces and related results by domains. Namely various theorems on traces concerning different domains will be collected in different sections.

The special attention in this paper will be given to recent research papers of the first author related with the trace problem and the plan of the paper is the following we provide first various definitions of various interesting analytic function spaces in various domains in ${\rm {\mathbb C}}^{n} $ then provide some precise information on traces of such type function spaces in ${\rm {\mathbb C}}.$ Note also that the trace operators may be also considered in relation with various embedding theorems and such type sharp results will also be given in this expository paper.

This is the first part of our expository paper on Traces in analytic function spaces of several complex variables which contains various new and old results on traces in the unit polydisk and in the unit polyball and some discussion of related problems.

\section{DIAGONAL MAPPINGS IN THE UNIT POLYDISK}

In this section we consider some interesting problems on diagonal mappings in the unit polydisk.

Let $U^{n} =\left\{z=(z_{1} ,\ldots ,z_{n} ):\left|z_{j} \right|<1,{\kern 1pt} j=\overline{1,n}\right\}$ be a the open unit polydisk in complex $n$-space ${\rm {\mathbb C}}^{n} $, $T^{n} $ the $n$-dimensional torus forming the distinguished boundary of $U^{n} $.

Let also $H(U^{n} )$ the space of functions holomorphic in $U^{n} $, $h^{p} (U^{n} )$ the Banach space of functions \textit{u }which are \textit{n}-harmonic (i.e. harmonic in each variable) in $U^{n} $ and which satisfy the condition 

\[\left\| u\right\| _{h^{p} } =\mathop{\sup }\limits_{0<r<1} \left(\int _{T^{n} } \left|u(r\zeta )\right|^{p} dm_{n} (\zeta )\right)^{\frac{1}{p} } <+\infty ,\] 

here $m_{n} $ is the normalized $n$-dimensional Lebesque measure on the torus $T^{n} $. We denote by $H^{p} (U^{n} )$ the closed subspace of $h^{p} (U^{n} )$ consisting of functions holomorphic in $U^{n} $.

Let $\left\{\mu _{j} \right\}_{j=1}^{n} $ be the finite Borel measures on $U,$ $\mu =(\mu _{1} ,...,\mu _{n} )$. We denote $L^{\vec{p}} (\mu )$ the space of $\mu $-measurable functions $h$ in $U^{n} $ for which 

\[\left\| h\right\| _{L^{\vec{p}} (\mu )} =\left(\int _{U} \left(\int _{U} ...\left(\int _{U} \left|h(z_{1} ,...,z_{n} )\right|^{p_{1} } d\mu _{1} \right)^{\frac{p_{2} }{p_{1} } } ...d\mu _{n-1} \right)^{\frac{p_{n} }{p_{n-1} } } d\mu _{n} \right)^{\frac{1}{p_{n} } } <+\infty .\] 

In paper [50] describe the positive bounded measures $\mu $ on the disk $U=U^{1} $ for which the diagonal operator $Df(z)=f(z,...,z),$ $z\in U,$ maps the space $H^{p} (U^{n} )$ into the usual $L^{\vec{p}} (\mu ),$ $0<p<+\infty .$ Allows one to give a complete description of the set $H^{p} (U^{n} )$ on the diagonal of the polydisk $U^{n} $.

This problem was formulated in Rudin's book [1]. Long history of the issue can be found in the review [9] and in the work of [2], [51].

For $0<p<+\infty ,$ $-1<\alpha <+\infty ,$ let $A_{\alpha }^{p} $ be the space of analytic functions on $U$ satisfying 

\[\left\| f\right\| _{A_{\alpha }^{p} } =\left(\iint \nolimits _{U} \left|f(z)\right|^{p} (1-\left|z\right|)^{\alpha } dxdy\right)^{\frac{1}{p} } <+\infty ,z=x+iy.\] 

W. Rudin proved that $\left\{Df:f\in H^{2} (U^{2} )\right\}=A_{0}^{2} $ and $\left\{Df:f\in H^{1} (U^{2} )\right\}=A_{0}^{1} $.

F.A. Shamoyan in [50] proved that $\left\{Df:f\in H^{p} (U^{n} )\right\}=A_{n-2}^{p} ,0<p<+\infty .$

Let $\zeta \in T^{n} ,$ $l\in {\rm {\mathbb R}}^{n} ,$ $l=(l_{1} ,...,l_{n} ),$ $0<l_{j} <1,$ $1\le j\le n,$ 

\[\Delta _{l} (\zeta )=\left\{z\in U^{n} :1-l_{j} \le \left|z_{j} \right|<1,\left|\arg z_{j} -\arg \zeta _{j} \right|\le \frac{l_{j} }{2} ,j=\overline{1,n}\right\}.\] 

\textbf{Theorem 1} ([50])\textbf{.} \textit{Let $f\in H^{p} ,$ $0<p\le q<+\infty ,$ $\mu $ be a probability measure on $U$. Then the following are equivalent:}

\textit{1. $\left(\int _{U} \left|Df\right|^{q} d\mu \right)^{\frac{1}{q} } \le c\left(\int _{T^{n} } \left|f\right|^{p} dm_{n} \right)^{\frac{1}{p} } ;$}

\begin{flushleft}
\textit{2. $\mu (\Delta _{l} (\zeta ))\le Cl^{\frac{nq}{p} } ,\zeta \in T^{1} ,0<l\le 1.$ }
\end{flushleft}

Note that it will be nice to obtain this result in bounded strongly pseudoconvex domains and in tubular domains over symmetric cones which is surely more difficult problem. This will be discussed in the second part of our expository paper.

Obviously, the measure $d\mu (r,\phi )=(1-r)^{n-2} rdrd\phi $ satisfies condition 2). Therefore $D$ maps $H^{p} (U^{n} )$ into $A_{n-2}^{p} ,$ $0<p<+\infty $.

The converse also holds, namely

\textbf{Theorem 2} ([50])\textbf{.} \textit{For each $f\in A_{n-2}^{p} ,$ $0<p<+\infty ,$ a function $g\in H^{p} (U^{n} )$ may be constructed such that $Dg(z)=f(z),$ $z\in U.$} 

Note that it will be nice to obtain this result in bounded strongly pseudoconvex domains and in tubular domains over symmetric cones which is surely more difficult problem.This will be discussed in the second part of our expository paper.

Combining theorems 1 and 2 we have the following subtle interesting result $Trace\, H^{p} =A_{n-2}^{p} $ for all positive $p,$ $n\ge 2.$

We also note that various partial results in this direction were obtained also by various experts A. Shields, P. Duren, D. Oberlin and Ch. Horowitz (we refer the interested reader for this for example to [2], and also to [9]).

Let $\omega _{j} $, $j=\overline{1,n},$ the non-negative functions of class $L^{1} (0,1),$ $0<p<+\infty .$

We denote $H^{p} (\omega _{1} ,...,\omega _{n} )=H^{p} (\vec{\omega })$ the class of holomorphic functions $f$ in $U^{n} $ for which 

\[\left\| f\right\| _{H^{p} (\omega _{1} ,...,\omega _{n} )} =\left(\int _{U^{n} } \left|f(z_{1} ,...,z_{n} )\right|^{p} \omega _{1} (1-\left|z_{1} \right|)...\omega _{n} (1-\left|z_{n} \right|)dm_{2n} (z_{1} ,...,z_{n} )\right)^{\frac{1}{p} } <+\infty ,\] 

$dm_{2n} $--$2n$-dimensional Lebesgue measure on $U^{n} .$

By the symbol $\Omega $ we denote the set of measurable non-negative on $\left(0,{\kern 1pt} 1\right)$ functions $\omega $, for which there exist positive numbers $m_{\omega } $, $M_{\omega } $, $q_{\omega } $, moreover $m_{\omega } $, $q_{\omega } $, such that $m_{\omega } \le \frac{\omega (\lambda r)}{\omega (r)} \le M_{\omega } $ for all $r\in (0,{\kern 1pt} {\kern 1pt} 1)$, $\lambda \in \left[q_{\omega } ,{\kern 1pt} 1\right]$.

Let $\omega _{j} \in \Omega ,$ $j=\overline{1,n},$ $0<p<+\infty ,$ $D(f)(z)=f(z,\ldots ,z),$$\Omega _{n} (r)=r^{2n-2} \prod \limits _{j=1}^{n} \omega _{j} (r),r\in (0,1).$

Let $H^{p} (\Omega _{n} )$ the class of holomorphic functions $f$ in $U$ for which 

\[\left\| f\right\| _{H^{p} (\Omega _{n} )} =\left(\int _{U} \left|f(z)\right|^{p} \Omega _{n} (1-\left|z\right|)dm_{2} (1-\left|z\right|)\right)^{\frac{1}{p} } <+\infty .\] 

\textbf{Theorem 3} ([8])\textbf{.} \textit{Let $\omega _{j} \in \Omega ,$ $j=\overline{1,n},$ $0<p<+\infty $. Then for any function $f\in H^{p} (\omega _{1} ,...,\omega _{n} ),$ $D(f)\in $ $H^{p} (\Omega _{n} )$, moreover $\left\| D(f)\right\| _{H^{p} (\Omega _{n} )} \le C(\omega _{1} ,...,\omega _{n} )\left\| f\right\| _{H^{p} (\omega _{1} ,...,\omega _{n} )} .$}

\textit{Conversely, for any $g\in H^{p} (\Omega _{n} ),$ we can construct $f\in H^{p} (\omega _{1} ,...,\omega _{n} )$ such that $D(f)(z)=g(z),$ $z\in U,$ that is $DH^{p} (\omega _{1} ,...,\omega _{n} )=H^{p} (\Omega _{n} ),$ $0<p<+\infty .$  }

Note that it will be nice to obtain this result in bounded strongly pseudoconvex domains and in tubular domains over symmetric cones which is surely more difficult problem. This will be discussed in the second part of our expository paper.

In the following theorem we provide complete descriptions of traces of general Bergman spaces (with general weights). In theorem 4 below even more general mixed norm Bergman space case is considered by us. Traces of Bergman spaces were studied also in J. Shapiro and some history can be seen also in [9].

Note also that very similarly trace problem can be posed in spaces of Bergman of pluriharmonic functions in the unit polydisk. And here also sharp results in this direction concerning traces can be also obtained as in classical analytic Bergman spaces in the unit polydisk. We leave this not difficult task to interested readers omitting details.

Let $L^{\vec{p}} (\vec{\omega })$ space of measurable functions $f$ in $U^{n} $ for which

\[\left\| f\right\| _{L^{\vec{p}} (\vec{\omega })} =(\int _{U} \omega _{n} (1-\left|\zeta _{n} \right|)...(\int _{U} \omega _{2} (1-\left|\zeta _{2} \right|)(\int _{U} \left|f(\zeta _{1} ,...,\zeta _{n} )\right|^{p_{1} } \times \] 

\[\times \omega _{1} (1-\left|\zeta _{1} \right|)dm_{2} (\zeta _{1} ))^{\frac{p_{2} }{p_{1} } } dm_{2} (\zeta _{2} ))^{\frac{p_{3} }{p_{2} } } ...dm_{2} (\zeta _{n} ))^{\frac{1}{p_{n} } } <+\infty ,\] 

where $U=U^{1} $, $dm_{2} (\zeta )$ -- 2-dimensional Lebesgue measure on $U$.

Let also $H^{\vec{p}} (\vec{\omega })=L^{\vec{p}} (\vec{\omega })\cap H(U^{n} ).$

It is not difficult to see that if $f\in H^{\vec{p}} (\vec{\omega })$ the function $\phi (z)=D(f)(z)=f(z,...,z),z\in U,$ is a holomorphic function in $U$. The problem of constructing a diagonal mapping in this space was solved in [38, 40].

Let $\vec{p}=(p_{1} ,...,p_{n} ),$ $\vec{\omega }=(\omega _{1} ,...,\omega _{n} ),$ $\tilde{\omega }_{p} (t)=\omega _{n} (t)\prod \limits _{j=1}^{n-1} \left(\omega _{j} (t)t^{2} \right)^{\frac{p_{n} }{p_{j} } } ,t\in (0,1].$

Let $H^{p_{n} } (\tilde{\omega }_{p} )$ the class of holomorphic functions $f$ in $U$ for which 

\[\left\| f\right\| _{H^{p_{n} } (\tilde{\omega }_{p} )} =\left(\int _{U} \left|f(z)\right|^{p_{n} } \tilde{\omega }_{p} (1-\left|z\right|)dm_{2} (z)\right)^{\frac{1}{p_{n} } } <+\infty .\] 

\textbf{Theorem 4} ([38])\textbf{.} \textit{Let $f\in H^{\vec{p}} (\vec{\omega }),$ $\vec{\omega }=(\omega _{1} ,...,\omega _{n} ),$ $\omega _{j} \in \Omega ,$ $\vec{p}=(p_{1} ,...,p_{n} ),$ $j=\overline{1,n}.$ Then $D(f)\in H^{p_{n} } (\tilde{\omega }_{p} ),$ moreover $\left\| D(f)\right\| _{H^{p_{n} } (\tilde{\omega }_{p} )} \le const\left\| f\right\| _{H^{\vec{p}} (\tilde{\omega })} .$}

\textit{Inversely, for any $f\in H^{p_{n} } (\tilde{\omega }_{p} )$, and $\alpha $ is a large enough positive number, }

\[F(z_{1} ,...,z_{n} )=\frac{\alpha n-1}{\pi } \int _{U} \frac{\left(1-\left|\zeta \right|^{2} \right)^{\alpha n-2} f(\zeta )dm_{2} (\zeta )}{(1-\overline{\zeta _{1} }z_{1} )^{\alpha } ...(1-\overline{\zeta _{n} }z_{n} )^{\alpha } } \] 

\textit{then $F\in H^{\vec{p}} (\tilde{\omega })$ that is $DH^{\vec{p}} (\tilde{\omega })=H^{p_{n} } (\tilde{\omega }_{p} ).$ }

Note that it will be nice to obtain this result in bounded strongly pseudoconvex domains and in tubular domains over symmetric cones which is surely more difficult problem. This will be discussed in the second part of our expository paper.

This theorem for $p_{j} =p,$ $j=\overline{1,n},$ was obtained in [8].

Let $L^{\vec{p}} (T^{n} )$,$\vec{p}=(p_{1} ,\ldots ,p_{n} )$, $0<p_{j} <+\infty $, $j=\overline{1,n},$ the space of measurable on $T^{n} $ functions $f$ for which the quasinorm is finite 

\[\left\| f\right\| _{L^{\vec{p}} (T^{n} )} =\left(\int _{T} \ldots \left(\int _{T} \left(\int _{T} \left|f(\zeta _{1} ,\ldots ,\zeta _{n} )\right|^{p_{1} } dm_{1} (\zeta _{1} )\right)^{\frac{p_{2} }{p_{1} } } dm_{1} (\zeta _{2} )\right)^{\frac{p_{3} }{p_{2} } } \cdots dm_{1} (\zeta _{n} )\right)^{\frac{1}{p_{n} } } ,\] 

where $T=T^{1} $, $dm_{1} $ the linear Lebesgue measure on $T$. 

 Let $\vec{p}=(p_{1} ,\ldots ,p_{n} )$,$0<p_{j} <+\infty $, $j=\overline{1,n}$, then we define generalized Hardy spaces with mixed norms $H^{\vec{p}} (U^{n} )$ as spaces of holomorphic in $U^{n} $ functions for which 

\[\mathop{{\rm s}up}\limits_{{\rm 0}\le {\rm r}<{\rm 1}} \; \left(\int _{T} \cdots \left(\int _{T} \left(\int _{T} \left|f(r\zeta )\right|^{p_{1} } dm_{1} (\zeta _{1} )\right)^{\frac{p_{2} }{p_{1} } } dm_{1} (\zeta _{2} )\right)^{\frac{p_{3} }{p_{2} } } \cdots dm_{1} (\zeta _{n} )\right)^{\frac{1}{p_{n} } } <\infty ,\] 

where $\zeta =(\zeta _{1} ,\ldots ,\zeta _{n} )\in T^{n} $. The case $p_{1} =\ldots =p_{n} =p$ space $H^{\vec{p}} (U^{n} )$ coincides with the usual Hardy space $H^{p} (U^{n} )$.

In [10-11] the characterization of traces of functions from spaces $H^{\vec{p}} (U^{2} )$ on the diagonal of the bidisk $U^{2} $ the authors obtained. That is, those functions $\phi \in H(U)$, for which the representation $\phi (z)=f(z,z)$, $f\in H^{p_{1} ,p_{2} } (U^{2} )$ is true.

Let $0<p,q<+\infty $, then 

\[A^{p,q} =\left\{f\in H(U),{\kern 1pt} {\kern 1pt} \left\| f\right\| _{A^{p,q} } =\left(\int _{0}^{1} \left(\int _{-\pi }^{\pi } \left|f(re^{i\phi } )\right|^{p} d\phi \right)^{\frac{q}{p} } dr\right)^{\frac{1}{q} } <+\infty \right\}.\] 

\textbf{Theorem 5} ([11])\textbf{.} \textit{Let $f\in H^{p_{1} ,p_{2} } ,$ $1\le p_{1} ,p_{2} <+\infty ,$ then function $\phi (z)=D(f)(z)=f(z,z)$ is of the class $A^{p_{1} ,p_{2} } $. On the contrary, if $\phi \in A^{p_{1} ,p_{2} } ,$ it is possible to construct explicitly the function $f\in H^{p_{1} ,p_{2} } (U^{2} )$ such that $\phi (z)=D(f)(z)=f(z,z),$ $z\in U,$ that is $D(H^{p_{1} ,p_{2} } )=A^{p_{1} ,p_{2} } $. }

Note that it will be nice to obtain this result in bounded strongly pseudoconvex domains and in tubular domains over symmetric cones which is surely more difficult problem. This will be discussed in the second part of our expository paper.

Some new interesting estimates on traces of new mixed norm analytic Hardy spaces but in the polyballs (products of unit balls) partially generalizing this interesting subtle result were obtained recently by the first author, but only for $p_{1} \le p_{2} \le p_{3} \le ...\le p_{n} $ in case based on properties of \textit{r}-lattices in the unit ball.

Note also that this trace problem related with mixed norm Hardy spaces in theorem 5 was solved only for $n=2$ particular case in simple unit polydisk.

The general case is an open interesting problem. Note also that similar problems can be posed also for very similar mixed norm Hardy spaces (analytic function spaces) but in context of very general tubular domains over symmetric cones and bounded strongly pseudoconvex domains with smooth boundary. We should note here that these trace problems surely are more complicated.

For $z=(z_{1} ,...,z_{n} )\in {\rm {\mathbb C}}^{n} ,$ $z_{j} =r_{j} e^{i\phi _{j} } ,$ $\zeta =(\zeta _{1} ,...,\zeta _{n} )\in {\rm {\mathbb C}}^{n} ,$ $\zeta _{j} =\rho _{j} e^{i\theta _{j} } ,$ $\alpha =(\alpha _{1} ,...,\alpha _{n} )\in {\rm {\mathbb R}}^{n} ,$ then $z^{\alpha } =z^{\alpha _{1} } \cdot ...\cdot z^{\alpha _{n} } ,$ $\left|\alpha \right|:=\alpha _{1} +...+\alpha _{n} ,$ $(1-\left|z\right|^{2} )^{\alpha } :=\prod \limits _{j=1}^{n} (1-\left|z_{j} \right|^{2} )^{\alpha _{j} } ,$ $(1-\zeta z)^{\alpha +2} :=\prod \limits _{j=1}^{n} (1-\zeta _{j} z_{j} )^{\alpha _{j} +2} ,$ and for $\omega =(\omega _{1} ,...,\omega _{n} ),\omega _{j} \in S,1\le j\le n,$ $\omega _{\Pi } (1-r):=\prod \limits _{j=1}^{n} \omega _{j} (1-r_{j} ),$ $\omega _{\Pi }^{s} (1-r):=\prod \limits _{j=1}^{n} \omega _{j}^{s} (1-r_{j} ),$ $\forall s\in {\rm {\mathbb R}},$ $r=(r_{1} ,...,r_{n} )\in Q_{n} =[0;1)^{n} .$

We denote $L_{\omega }^{p,q} (U^{n} )$ the space of measurable functions $f$ in $U^{n} $ for which 

\[\left\| f\right\| _{L_{\omega }^{p,q} } =\left(\int _{Q^{n} } \omega _{\Pi } (1-r)\left(\int _{T^{n} } \left|f(r\omega )\right|^{p} dm_{n} (\omega )\right)^{\frac{q}{p} } dr\right)^{\frac{1}{q} } <+\infty ,0<p,q<+\infty ,\] 

where $dm_{n} $ is Lebesgue measure on $T^{n} $.

Let also $A_{\omega }^{p,q} (U^{n} )=H(U^{n} )\cap L_{\omega }^{p,q} (U^{n} )$.

Obviously, if $f\in H(U^{n} ),$ $f\in A_{\omega }^{p,q} (U^{n} )$, then $D(f)(z)=f(z,...,z)\in H(U).$

Let also $\Omega _{n} (r)=r^{(\frac{q}{p} +1)(n-1)} \prod \limits _{j=1}^{n} \omega _{j} (r),$ $r\in (0,1).$

Let $A_{\Omega _{n} }^{p,q} (U)$ the class of holomorphic functions $f$ in $U$ for which 

\[\left\| f\right\| _{A_{\Omega _{n} }^{p,q} } =\left(\int _{0}^{1} \Omega _{n} (1-r)\left(\int _{-\pi }^{\pi } \left|f(r\omega )\right|^{p} dm(\omega )\right)^{\frac{q}{p} } dr\right)^{\frac{1}{q} } <+\infty ,0<p,q<+\infty .\] 

In [14, 37] was obtained a complete characterization of the traces of weighted anisotropic spaces of analytic functions in the polydisk with mixed norm on the on the diagonal of the polydisk.

\textbf{Theorem 6} ([37])\textbf{.} \textit{Let $\omega =(\omega _{1} ,...,\omega _{n} ),$} $\omega _{j} \in \Omega ,$\textit{ $j=\overline{1,n},$ $0<p,q<+\infty .$ Then we have the following equality $Trace\, A_{\omega }^{p,q} (U^{n} )=A_{\Omega _{n} }^{p,q} (U).$}  

Note that it will be nice to obtain this result in bounded strongly pseudoconvex domains and in tubular domains over symmetric cones which is surely more difficult problem. This will be discussed in the second part of our expository paper.

This theorem was proved in less general case for less general weights in a paper of [52] using completely other technique.

The natural interesting nice open question here is to find complete descriptions of traces of such type analytic mixed norm spaces but in context of polyballs (products of balls) or tubular domains over symmetric cones or in pseudoconvex domains with smooth boundary.

First note that we can easily define such type analytic function spaces in mentioned domains and product of such type domains using basic definitions of analytic function spaces in such domains. Then to try to solve Trace problem in such type function spaces.

Some new interesting sharp results on BMOA type analytic function spaces and Herz type analytic function spaces on diagonal map, and various sharp theorems or extensions of known results concerning diagonal map in the polydisk can be seen in papers of first author [12, 15-18, 26]. The recent research paper of G. Ren and J. Shi [52] contains also sharp results on diagonal map in mixed norm function spaces. We refer the reader to these papers.

Note in addition that if $a_{k} $ is a $r$-lattice in the unit polydisk and $D(a_{k} ,r)$ is a Bergman ball there, then we can define new Herz type analytic spaces in the unit polydisk with quazinorms $\sum \limits _{k} (\int _{D(a_{k} ,r)} |f(z)|^{q} (1-|z|)^{\alpha } dm_{2n} (z))^{p} $ and pose a problem to find complete description of traces of these analytic spaces.

\section{TRACES IN VARIOUS SPACES OF ANALYTIC FUNCTIONS IN POLYBALL}

In this section we consider some interesting problems on traces in various spaces of analytic functions in polyball.

If we formally simply replace the unit disk in our previous section by the unit ball then we naturally arrive at the same problem as diagonal map in polydisk from previous section but here the task is to find precise estimates of traces for functions from analytic spaces in products of the unit balls (polyballs). And this problem is surely more general than what we had in previous section in the unit polydisk.

We collect in this section several such type sharp results concerning traces of various analytic function spaces in polyballs (mixed norm spaces, Bloch type, Bergman type, BMOA type spaces).

Let the open unit ball in ${\rm {\mathbb C}}^{n} $ is the set $B=B^{n} =\{ z\in {\rm {\mathbb C}}^{n} :\left|z\right|<1\} .$ The boundary of $B$ will be denoted by $S,$ $S=\{ z\in {\rm {\mathbb C}}^{n} :\left|z\right|=1\} $ (see, for example, [53]).

Let $d\upsilon $ denote the volume measure on $B$, normalized so that $\upsilon (B)=1$, and let $d\sigma $ denote the surface measure on $S,$\textit{ }normalized so that $\sigma (S)=1$.

For $\alpha >-1$ the weighted Lebesgue measure $d\upsilon _{\alpha } $ is defined by $d\upsilon _{\alpha } =c_{\alpha } \left(1-\left|z\right|^{2} \right)^{\alpha } d\upsilon (z),$ where $c_{\alpha } =\frac{\Gamma (n+\alpha +1)}{\Gamma (n+1)\Gamma (\alpha +1)} .$ Is a normalizing constant, so that $\upsilon _{\alpha } (B)=1$ (see [54]).

Let also $d\tilde{\upsilon }_{\alpha } (z)=d\upsilon _{\alpha } (z_{1} )...d\upsilon _{\alpha } (z_{m} )=\left(1-\left|z_{1} \right|^{2} \right)^{\alpha } ...\left(1-\left|z_{m} \right|^{2} \right)^{\alpha } d\upsilon (z_{1} )...d\upsilon (z_{m} )$.

Let every functions $f\in H(B)$ having a series expansion $f(z)=\sum \limits _{\left|k\right|\ge 0} a_{k} z^{k} ,$ we define the operator of fractional differentiation by ${\rm {\mathcal D}}^{\alpha } f(z)=\sum \limits _{\left|k\right|\ge 0} \left(\left|k\right|+1\right)^{\alpha } a_{k} z^{k} ,$ where $\alpha $ is any real number.

For $z\in B$ and $r>0$ the set $D(z,r)=\left\{\omega \in B:\rho (z,\omega )<r\right\},$ where $\rho $ is a Bergman metric on $B$, $\rho (z,\omega )=\frac{1}{2} \log \frac{1+\left|\phi _{z} (\omega )\right|}{1-\left|\phi _{z} (\omega )\right|} $ is called the Bergman metric ball at $z$ (see [54]), where $\phi _{z} $ is the M$\ddot{o}$bius transformation of $B$ that interchanges 0 and $z$.

For $\alpha >-1$ and $p>0$ the weighted Bergman space $A_{\alpha }^{p} $ consists of holomorphic functions f in $L^{p} (B,d\upsilon _{\alpha } )$ that is, $A_{\alpha }^{p} (B)=L^{p} (B,d\upsilon _{\alpha } )\cap H(B).$

Let $B^{m} =B_{n}^{m} =B\times ...\times B$ be a polyball (product of m balls), and let $H(B^{m} )$ be a space of all holomorphic functions $f(z_{1} ,...,z_{m} ),$ $z_{j} \in B,$\textit{$j=\overline{1,m}.$}

Bergman classes on polyballs $A^{p} (B^{m} ,d\upsilon _{\alpha _{1} } ...d\upsilon _{\alpha _{m} } )$ consists of function $f$ in $H(B^{m} )$, such that 

\[\int _{B} ...\int _{B} \left|f(z_{1} ,...,z_{m} )\right|^{p} (1-\left|z_{1} \right|)^{\alpha _{1} } ...(1-\left|z_{m} \right|)^{\alpha _{m} } d\upsilon (z_{1} )...d\upsilon (z_{m} )<+\infty ,\] 

\[0<p<+\infty , \alpha _{j} >-1, j=\overline{1,m},\] 

As in previous section very similarly we can pose a trace problem but not in the unit polydisk but in the unit polyball.

We observe that for $n=1$ this problem completely coincide with the well-known problem of diagonal map. The last problem of description of diagonal of various subspaces of $H(U^{n} )$ of spaces of all holomorphic functions in the polydisk was studied by many authors before (see, for example, [1-2, 9, 52] and references there).

We formulate below new sharp trace theorems for Bergman and Bloch type analytic function spaces in the polyball.

\textbf{Theorem 7} ([24])\textbf{.} \textit{Suppose $0<p\le +\infty ,$ $s_{1} ,...,s_{m} >-1.$ Put $t=(m-1)(n+1)+\sum \limits _{j=1}^{m} s_{j} .$ Then $Trace\left(A^{p} (B^{m} ,d\upsilon _{s_{1} } ...d\upsilon _{s_{m} } )\right)=$ $A^{p} (B,d\upsilon _{t} ).$} 

For any positive real numbers $r_{j} $, \textit{$j=\overline{1,m},$} we define 

\[\Lambda \left(r_{1} ,...,r_{m} \right):=\left\{f\in H(B^{m} ):\mathop{\sup }\limits_{z_{j} \in B} \left(\left|f(z_{1} ,...,z_{m} )\right|\prod \limits _{j=1}^{m} \left(1-\left|z_{j} \right|^{2} \right)^{r_{j} } \right)<+\infty \right\}.\] 

\textbf{Theorem 8} ([24])\textbf{.} \textit{Let $r_{j} >0,$ $j=\overline{1,m},$ $r=r_{1} +...+r_{m} $. Then $Trace\left(\Lambda (r_{1} ,...,r_{m} )\right)=\Lambda (r).$ }

Let 

\[\Lambda _{\log } \left(r_{1} ,...,r_{m} \right):=\left\{f\in H(B^{m} ):\mathop{\sup }\limits_{z_{j} \in B} \left(\left|f(z_{1} ,...,z_{m} )\right|\prod \limits _{j=1}^{m} \left(\log \frac{1}{1-\left|z_{j} \right|} \right)^{-\frac{1}{r_{j} } } \left(1-\left|z_{j} \right|\right)^{\frac{1}{r_{j} } } \right)<+\infty \right\},\] 

\[\sum \limits _{j=1}^{m} \frac{1}{r_{j} } =1, r_{j} >0,j=\overline{1,m}.\] 

\[\Lambda _{\log } \left(1\right):=\left\{f\in H(B):\mathop{\sup }\limits_{z\in B} \left(\left|f(z)\right|\left(\log \frac{1}{1-\left|z\right|} \right)^{-1} \left(1-\left|z\right|\right)\right)<+\infty \right\}.\] 

\textbf{Theorem 9} ([24])\textbf{.} \textit{$Trace\left(\Lambda _{\log } (r_{1} ,...,r_{m} )\right)=\Lambda _{\log } \eqref{GrindEQ__1_}.$  }

Let 

\[A_{\alpha _{1} ,...,\alpha _{m} }^{p_{1} ,...,p_{m} } \left(B^{m} \right)=\left\{f\in H(B^{m} ):\right. \left\| f\right\| _{A_{\alpha _{1} ,...,\alpha _{m} }^{p_{1} ,...,p_{m} } } :=\] 

\[:=\left(\int _{B} \left(1-\left|z_{m} \right|\right)^{\alpha _{m} } \left(...\int _{B} \left|f\left(z_{1} ,...,z_{m} \right)\right|^{p_{1} } \left(1-\left|z_{1} \right|\right)^{\alpha _{1} } d\upsilon \left(z_{1} \right)\right)^{\frac{p_{2} }{p_{1} } } ...d\upsilon \left(z_{m} \right)\right)^{\frac{1}{p_{m} } } <+\infty ,\] 

\[0<p_{j} <+\infty , \alpha _{j} >-1, j=\overline{1,m}.\] 

\textbf{Theorem 10 }([28])\textbf{.} \textit{Let $\gamma =\alpha _{m} +\sum \limits _{j=1}^{m-1} \left(n+1+\alpha _{j} \right)\frac{p_{m} }{p_{j} } ,$ $\alpha _{j} >-1,$ $p_{j} >1,$ $j=\overline{1,m}.$ Then $Trace\left(A_{\alpha _{1} ,...,\alpha _{m} }^{p_{1} ,...,p_{m} } \left(B^{m} \right)\right)=A_{\gamma }^{p_{m} } \left(B\right).$ }

Let ${\rm {\mathbb C}}_{+} $ be the upper half-plane and ${\rm {\mathbb C}}_{+}^{m} ={\rm {\mathbb C}}_{+} \times ...\times {\rm {\mathbb C}}_{+} $ (\textit{m}-times).

We provide below a new sharp trace theorem for Bergman spaces in products of upper half spaces.

\textbf{Theorem 11 }([28])\textbf{.} \textit{Let $0<p\le +\infty $ and} $\alpha _{1} ,...,\alpha _{m} >-1,$\textit{ $\gamma =\sum \limits _{j=1}^{m} \alpha _{j} +2m-2.$ Then $Trace(A_{\alpha _{1} ,...,\alpha _{m} }^{p} ({\rm {\mathbb C}}_{+}^{m} ))=A_{\gamma }^{p} ({\rm {\mathbb C}}_{+} ).$} 

We define now below new analytic BMOA type spaces in the polyball and provide a new sharp trace theorem for analytic BMOA type spaces in the polyball. Later this result was extended to more difficult tube and bounded strongly pseudoconvex domain by first author based on arguments provided in the simpler case of polyball.

This sharp theorem can be also formulated and proved for rather general mixed norm Bergman classes $A_{\alpha _{1} ,...,\alpha _{m} }^{p_{1} ,...p_{m} } $ in product of upper spaces ${\rm {\mathbb C}}_{+}^{m} $ very similarly as in the unit polyball case which was formulated by us above.

Let $0<p<+\infty .$ For the real numbers $s_{1} ,...,s_{m} >-1$ and $r_{1} ,...,r_{m} >0$ we define $M_{r_{1} ,...,r_{k} }^{p} (B_{n}^{m} ,d\upsilon _{s_{1} } ...d\upsilon _{s_{m} } )$ to be space of all measurable functions $f$ on $B_{n}^{m} $ for which the measure $d\mu _{f} =\left|f\right|^{p} d\upsilon _{s_{1} } ...d\upsilon _{s_{m} } $ is a bounded $(r_{1} ,...,r_{m} )$-Carleson measure.

For any $f\in M_{r_{1} ,...,r_{k} }^{p} (B_{n}^{m} ,d\upsilon _{s_{1} } ...d\upsilon _{s_{m} } )$ we define 

\[\left|\left\| f\right\| \right|^{p} =\sup \left(\frac{\mu _{f} \left(E(a_{1} \right)\times ...\times E(a_{m} ))}{\left(1-\left|a_{1} \right|\right)^{r_{1} } ...\left(1-\left|a_{m} \right|\right)^{r_{m} } } :a_{1} ,...,a_{m} \in B_{n} \right)<+\infty .\] 

Let $HM_{r_{1} ,...,r_{k} }^{p} (B_{n}^{m} ,d\upsilon _{s_{1} } ...d\upsilon _{s_{m} } )=H(B_{n}^{m} )\cap M_{r_{1} ,...,r_{k} }^{p} (B_{n}^{m} ,d\upsilon _{s_{1} } ...d\upsilon _{s_{m} } )$.

\textbf{Theorem 12} ([29])\textbf{.} \textit{Let $0<p<+\infty .$ Suppose that $s_{1} ,...,s_{m} >-1$ and $r_{1} ,...,r_{m} >0$ are such that $r_{j} <n+1+s_{j} $ for $j=\overline{1,m}.$ Put $t=(m-1)(n+1)+\sum \limits _{j=1}^{m} s_{j} $ and $r=\sum \limits _{j=1}^{m} r_{j} .$}

\textit{Then $Trace\left(HM_{r_{1} ,...,r_{k} }^{p} (B_{n}^{m} ,d\upsilon _{s_{1} } ...d\upsilon _{s_{m} } )\right)=HM_{r}^{p} (B_{n} ,d\upsilon _{t} ).$ }

In [22-24] extend results on trace problem for functions holomorphic on polyballs and give descriptions of traces for several concrete functional classes on polyballs defined with the help of Bergman metric ball.

We now formulate below another sharp result the complete analogue of sharp trace theorem~10 but for $p_{j} \le 1,$ \textit{$j=\overline{1,m},$} case.

 \textbf{Theorem 13} ([22-24])\textbf{.} \textit{Let $\gamma =\alpha _{m} +\sum \limits _{j=1}^{m-1} \left(n+1+\alpha _{j} \right)\frac{p_{m} }{p_{j} } ,$ $\alpha _{j} >-1,$ $0<p_{j} \le 1,$$j=\overline{1,m}.$ Then $Trace\left(A_{\alpha _{1} ,...,\alpha _{m} }^{p_{1} ,...,p_{m} } \left(B^{m} \right)\right)=A_{\gamma }^{p_{m} } \left(B\right).$}

Let $\alpha _{j} >-1,$ $\beta _{j} >-1,$ \textit{$j=\overline{1,m},$}$\beta >-1,$ $0<p,q<+\infty .$

In [22-24] define Herz type spaces of analytic functions in polyballs 

\[\left(K_{\alpha ,\beta }^{p,q} \right)=\left\{f\in H(B^{m} ):\int _{B} ...\int _{B} \left(\int _{D(z_{1} ,r)} ...\int _{D(z_{m} ,r)} \left|f(z)\right|^{p} d\tilde{\upsilon }_{\alpha } (z)\right)^{\frac{q}{p} } d\tilde{\upsilon }_{\beta } (z)<+\infty \right\};\] 

\[\left(\tilde{K}_{\alpha ,\beta }^{p,q} \right)=\left\{f\in H(B^{m} ):\int _{B} \left(\int _{D(z,r)} ...\int _{D(z,r)} \left|f(z)\right|^{p} d\tilde{\upsilon }_{\alpha } (z)\right)^{\frac{q}{p} } d\upsilon _{\beta } (z)<+\infty \right\}.\] 

We give below new sharp trace theorems concerning new analytic Herz type spaces in the polyball. These results can be probably expanded to more difficult domains.

\textbf{Theorem 14} ([22-24])\textbf{.} \textit{Let $0<p\le q\le 1,$ $t_{j} >-1,$ $\beta _{j} >-1,$ $j=\overline{1,m},$ $\alpha >-1,$ and }

\begin{flushleft}
\textit{$\alpha =\sum \limits _{j=1}^{m} \left[\left(n+1+\beta _{j} \right)\frac{p}{q} +(n+1+t_{j} )\right]-(n+1).$ Then $Trace\left(K_{\beta ,t}^{q,p} \left(B^{m} \right)\right)=A_{\alpha }^{p} \left(B\right).$}
\end{flushleft}

\textbf{Theorem 15} ([22-24])\textbf{.} \textit{Let $q>1,$ $p\le q,$ $t_{j} >-1,$ $\beta _{j} >-1,$ $j=\overline{1,m},$ $\alpha >-1,$ and }

\begin{flushleft}
\textit{$\alpha =\sum \limits _{j=1}^{m} \left[\left(n+1+\beta _{j} \right)\frac{p}{q} +(n+1+t_{j} )\right]-(n+1).$ Then $Trace\left(K_{\beta ,t}^{q,p} \left(B^{m} \right)\right)=A_{\alpha }^{p} \left(B\right).$}
\end{flushleft}

\textbf{Theorem 16} ([22-24])\textbf{.} \textit{Let $q>1$ or $q\le 1,$ $p\le q,$ $t>-1,$ $\beta _{j} >-1,$ $j=\overline{1,m},$ $\alpha >-1,$ and $\alpha =\sum \limits _{j=1}^{m} \left[\left(n+1+\beta _{j} \right)\frac{p}{q} \right]+t.$ Then $Trace\left(\tilde{K}_{\beta ,t}^{q,p} \left(B^{m} \right)\right)=A_{\alpha }^{p} \left(B\right).$}

Theorems 14, 15 and 16 for $p=q$ were proved in [23-24, 27-28].

Note that we defined a fractional derivative of analytic function in the ball in the start of this section and we can define similarly the same type operator of fractional derivative also for analytic functions on products of the unit balls. This task is very easy.

As a result we can also define for various $X$ analytic function spaces on polyballs more general the $D^{\alpha } $$X$ analytic function spaces in polyballs for example for positive meaning analytic function spaces consisting of analytic functions for which $D^{\alpha } f$ function belongs to$X$ space in the unit polyball. And we note that many our sharp results on traces in the unit polyballs of this section may be extended to such type analytic function classes in the polyballs also.

Some sharp results of this third section concerning polyballs can be via rather similar proofs extended partially to various more complicated analytic function spaces in tubular domains over symmetric cones and bounded strongly pseudoconvex domains with smooth boundary without difficulties.

Some descriptions (or estimates) of traces of analytic spaces $\left(M_{\alpha }^{p} \right),$ $\left(K_{\alpha ,\beta }^{p,q} \right),$ $\left(D_{\alpha ,\beta }^{p} \right)$ and related spaces 

\[\left(M_{\alpha }^{p} \right)=\left\{f\in H(B^{m} ):\int _{S} \int _{\Gamma _{t} (\zeta )} ...\int _{\Gamma _{t} (\zeta )} \left|f(z)\right|^{p} d\nu _{\alpha } (z)d\sigma (\zeta )<+\infty \right\},\] 

\[0<p<+\infty , \alpha >-1;\] 

\[\left(K_{\alpha ,\beta }^{p,q} \right)=\left\{f\in H(B^{m} ):\int _{B} ...\int _{B} \left(\int _{D(z_{1} ,r)} ...\int _{D(z_{m} ,r)} \left|f(z)\right|^{p} d\nu _{\alpha } (z)\right)^{\frac{q}{p} } d\nu _{\beta } (z)<+\infty \right\},\] 

\[0<p,q<+\infty , \alpha >-1, \beta >-1;\] 

\[\left(D_{\alpha ,\beta }^{p} \right)=\left\{f\in H(B^{m} ):\int _{0}^{1} (1-r)^{\beta } \left(\int _{\left|z_{1} \right|<r} ...\int _{\left|z_{m} \right|<r} \left|f(z)\right|^{p} \prod \limits _{j=1}^{m} \left(1-\left|z_{j} \right|\right)^{\alpha _{j} } d\nu (z_{j} )\right)dr<+\infty \right\},\] 

\[0<p<+\infty , \beta >-1, \alpha _{j} >-1,j=\overline{1,m}.\] 

Can be seen in the papers [12, 23] we note first they were obtained in the polydisk by M.~Jevtic, M. Pavlovic and R. Shamoyan then later these sharp results were extended to polyballs by O. Mihic and R. Shamoyan.

Many related and (or) very close to these trace theorems interesting results and interesting problems were considered also by various authors in recent decades in various analytic spaces in various domains in ${\rm {\mathbb C}}^{n} ,$ we mention for example papers [55-59].

We will not discuss these problems or these type (trace type) problems in this expository paper referring the authors to these mentioned interesting papers.

We included for conveniences of our readers in the list of the references below various other interesting new and old results on traces of analytic function spaces in more complicated tubular domains over symmetric cones and bounded strongly pseudoconvex domains. These results will be discussed and included in the second part of our expository paper.

Note in addition that if $a_{k} $ is a $r$-lattice in the unit polyball and $D(a_{k} ,r)$ is a Bergman ball there, then we can define new Herz type analytic spaces in the unit polyball with quazinorms $\sum \limits _{k} (\int _{D(a_{k} ,r)} |f(z)|^{q} (1-|z|)^{\alpha } dm_{2n} (z))^{p} $ and pose a problem to find complete description of traces of these analytic spaces.

Some very interesting extensions of diagonal map (trace) problem in the unit polydisk based on finite Blaschke products in the unit polydisk was considered by D. Clark [49]. Some very interesting applications of diagonal map can be seen in a paper of E. Amar and C. Menini [48].

We mention separately also very interesting paper [25] where Bergman and Bloch type analytic spaces on expanded disk and subframe were considered and some diagonal map type maps were studied.

We also wish to note that in recent monographs [2, 54, 59, 61] many interesting new (and not new) analytic function spaces were studied, they can be extended to product domains in a natural way and trace problem for such type analytic spaces in product domains also can be posed.

In [62-63] new interesting trace theorems and related to them sharp embedding theorems were considered in spaces of harmonic functions, this will be discussed in the second part of this paper.

Note also, some interesting technique from [41] can be used to get sharp (or not) trace theorems in various analytic function spaces in bounded symmetric domains.

Some interesting relations of trace theorems with Martinelly-Bochner integrals can be seen in a paper of the first author and S. Kurilenko [19]. Note also that several results of this paper later completely were extended to more general tube and pseudoconvex domains, these papers are also included in references by us.

\section{CONCLUSIONS}

Trace theorems in analytic function spaces in the polydisk and polyball  may have various interesting applications in complex function theory of several complex variables and also probably can be applied to solutions of various interesting problems in operator theory related to function spaces in product domains.

This article is the first part of our expository paper, the second part of this note will contain in particular some direct generalizations of results of this paper to more complicated domains as tubular domains over symmetric cones and bounded pseudoconvex domains.

\end{document}